\documentclass[preprint,review,10pt]{elsarticle}

\usepackage{graphicx}
\usepackage{natbib}
\usepackage{psfrag}
\usepackage{subfigure}
\usepackage{url}
\usepackage{epsfig}
\usepackage{amsmath}
\usepackage{algorithm}
\usepackage{algorithmic}
\usepackage{amssymb}
\usepackage{amsthm}
\usepackage{fancybox}

\newcommand{\tb}{\textbf{t}}
\newcommand{\vb}{\textbf{v}}
\newcommand{\xb}{\textbf{x}}
\newcommand{\ub}{\textbf{u}}
\newcommand{\ld}{\lambda}

\newcommand{\rset}{\mathbb{R}}

\newcommand{\bmtrx}{\left[\begin{array}}
\newcommand{\emtrx}{\end{array}\right]}

\newcommand{\Pp}{\mathbb{P}}

\newtheorem{theorem}{Theorem}[section]

\newtheorem{assumption}[theorem]{Assumption}

\newenvironment{remark}{
       \refstepcounter{theorem}\begin{trivlist}\item[]{\bfseries
       Remark \thetheorem\,}}
       {\end{trivlist}}

\newenvironment{example}{
       \refstepcounter{theorem}\begin{trivlist}\item[]{\bfseries
       Example \thetheorem\,}}
       {\end{trivlist}}

\begin{document}
\begin{frontmatter}

\title{ Parallel and distributed optimization methods for estimation and control in networks \tnoteref{acknow}}

\tnotetext[acknow]{ The research leading to these results has
received funding from: the European Union, Seventh Framework
Programme (FP7/2007--2013) under grant agreement no 248940;
CNCSIS-UEFISCSU (project TE, no. 19/11.08.2010); ANCS (project PN
II, no. 80EU/2010); Sectoral Operational Programme Human Resources
Development 2007-2013 of the Romanian Ministry of Labor, Family and
Social Protection through the Financial Agreement
POSDRU/89/1.5/S/62557.}

\author[First]{Ion Necoara}
\author[First]{Valentin Nedelcu}
\author[First]{Ioan Dumitrache}

\address[First]{University Politehnica Bucharest, Automatic Control and Systems Engineering
Department, 060042 Bucharest, Romania\\
(e-mail: ion.necoara,ioan.dumitrache,v.nedelcu@upb.ro)}

\begin{abstract}
System performance  for networks composed of interconnected
subsystems can be increased if the traditionally  separated
subsystems are jointly optimized. Recently, parallel and distributed
optimization methods have emerged as a powerful tool for solving
estimation and control problems in large-scale networked systems. In
this paper we review and analyze the optimization-theoretic concepts
of  parallel and distributed  methods for solving coupled
optimization problems and demonstrate how several  estimation and
control problems related to complex networked systems can be
formulated in these settings. The paper presents a systematic
framework for exploiting  the potential of the  decomposition
structures as a way to obtain different  parallel algorithms, each
with a different tradeoff among convergence speed, message passing
amount and distributed computation architecture. Several specific
applications from estimation and process control are included to
demonstrate the power of the approach.
\end{abstract}


\end{frontmatter}

%%%%%%%%%%%%%%%    INTRODUCTION     %%%%%%%%%%%%%

\section{Introduction}
\label{intro} In many application fields, the notion of networks has
emerged as a central, unifying concept for solving different
problems  in systems and control theory such as analysis, process
control and estimation.  We live and operate in a networked world.
We drive to work on  networks of roads and communicate with each
other using an  elaborate set of devices such as phones or
computers, that connect wirelessly and through the internet.
Traditional networks include transportation networks (roads, rails)
and networks of utilities (water, electricity, gas). But more recent
examples of the increasing impact of networks include information
technology networks (internet, mobile phones, acoustic networks,
etc), information networks (co-author networks, bibliographic
networks), social networks (collaborations, organizations), and
biological and genetic networks.

These networks are often composed of multiple subsystems
characterized by complex dynamics and mutual interactions such that
local decisions have long-range effects throughout the entire
network. Many problems associated to  networked systems, such as
state estimation and  control, can be posed as \textit{coupled
optimization problems} (see e.g.
\cite{CamJia:02,EbrBal:00,FarFer:10,KevBor:06,MaePen:10,
NegSch:08,NecSuy:08, VenRaw:05}, etc). Note that in these systems
the interaction between subsystems gives rise to coupling in the
cost or constraints, but with a specific algebraic structure, in
particular sparse matrix representation that could be exploited in
numerical algorithms. Therefore, in order to design an overall
decision architecture for such complex networks we need to solve
large coupled optimization problems but with specific structure. The
major difficulty in  these problems is that due to their size,
communication restrictions, or requirements on robustness, often no
central decisions can be taken; instead, the decisions have to be
taken locally. In such a set-up, single units, or local agents, must
solve local optimization subproblems and then they must negotiate
their outcomes and requirements with their neighbors in order to
achieve convergence to the global optimal solution. Basically, there
are two general
optimization approaches:\\
(i)  ``Centralized'' optimization algorithms: In this class the
specific structure of the system  is exploited, as it represents
considerable sparsity in the optimization problem due to the local
coupling between optimization variables (sometimes referred to as
\textit{separable optimization problems}).  The sparsity of the
problem, given by the influences between the subsystems, leads to
coupling constraints represented by sparse matrices. Though parts of
the algorithms will be parallelized, the parallelization in these
algorithms is not restricted by e.g. limited communication between
subsystems and is just for the sake of exploiting sparsity.  In
summary, ``centralized'' algorithms benefit from the sparsity
induced by the networked system  and solve the resulting
optimization problems on a parallel computer architecture.  Several
standard parallel and distributed optimization methods can be found
in the textbooks \cite{BerTsi:89, Las:70,MesMac:70}. Various survey
papers also exist on optimization-based distributed control. In the
70's Tamura \cite{Tam:73} and Mahmoud  \cite{Mah:77} presented very
comprehensive overviews. More recently, in \cite{RawSte:08} the
actual status of research in the field of coordinated
optimization-based control is presented. Many different control
topologies can be considered in distributed control, which have been
reviewed recently in \cite{Sca:09}. When there is no need to solve
the separable optimization problem  on a parallel computer
architecture, an alternative would be to solve the global
optimization problem using sparse solvers that take into account the
sparse structure of the problem at the linear algebra level of the
optimization algorithm. In general, this choice could lead to faster
algorithms in terms of
 CPU time than distributed or parallel algorithms. \\
(ii)   Distributed optimization algorithms (sometimes referred to as
\textit{distributed multi-agent optimization algorithms}):  In
contrast to the  ``centralized'' algorithms, distributed algorithms
on graphs have to satisfy an extra constraint, namely their
computations shall be performed on all nodes in parallel, and the
communication between nodes is restricted to the edges of the graph,
i.e. such algorithms do not use all-to-all communication protocols.
In many complex networked systems the desired behavior can be
formulated as coupled optimization problems but with restrictions on
communication due to the special network topology: e.g.
  estimation in sensor networks, consensus and rendezvous problems
   in multi-agent systems, resource allocation in computer networks
   \cite{EbrBal:00,JohSpe:08,OlfFax:07}.
    Some existing distributed methods that take into account
     explicitly information restrictions in the network combine
      consensus negotiations (as an efficient method  for information fusion)
      with subgradient methods \cite{JohKev:08,NecDum:10,NedOzd:09,NedOzd:10,ZhuMar:10}.

The goal of this paper is twofold: (i)  to establish a relationship
between estimation and control  in networked systems and distributed
optimization methods and demonstrate  the effectiveness of utilizing
optimization-theoretic approaches for  controlling such complex
systems; (ii)  motivated by this connection, to build upon
optimization based results  to better accommodate    a broader class
of estimation and control problems.  The core of this paper consists
of Section \ref{sec_application}, covering three applications of
estimation and control   that appear in the context of networked
systems and then proving how we can reformulate them as coupled
optimization problems. One of the key contributions of this paper is
to provide an accessible, yet relatively comprehensive, overview of
three classes of decomposition schemes from mathematical programming
for solving distributively  coupled optimization problems. We
demonstrate how the decomposition schemes suggest network
architectures and protocols with different properties in terms of
convergence speed and coordination overhead.  We also present new
decomposition methods that are more efficient in terms of
convergence speed  than some classical decomposition schemes.

The paper is organized as follows. In Section \ref{sec_application}
we introduce  different  estimation and control problems   that
appear in the context of  complex systems with interacting
subsystems dynamics and then we show how we can reformulate them as
coupled optimization problems. In Section \ref{sec_decomposition} we
present several parallel and distributed methods for solving this
type of structured optimization problems and analyze their
performance. Section  \ref{sec_decomposition} thus serves both as a
review of the necessary background and a summary of our new
extensions on decomposition methods.  For each of the applications,
numerical experiments on different parallel and distributed
algorithms are provided.

%%%%%%%%%%%%%%%%   APPLICATIONS      %%%%%%%%%%

\section{Estimation and control problems in networks}
\label{sec_application}

In this section we formulate different estimation and control
problems for systems consisting of interconnected subsystems. In
Subsection \ref{subsec_estimation} we present a state estimation
problem for a system, using a network of sensors which must exchange
information  in order to reach a consensus on the state estimated
for the entire  system. In Subsection \ref{subsec_optimal} we will
present the problem of optimal control for a large-scale system,
whose subsystems are coupled with their neighbors but the objective
function is decoupled.  Finally, in Subsection
\ref{subsec_cooperative} and \ref{subsec_cooperative_coupled} we
will discuss the cooperative control problem for a group of systems
(agents), which have decoupled or coupled dynamics but share a
common goal.

%%%%%%%%%%%%%%   STATE    ESTIMATION  %%%%%%%%%%%%%

\subsection{State estimation problem}
\label{subsec_estimation}
 In this section we formulate the distributed
state estimation problem for systems using a  sensor network  based
on the  moving horizon estimation (MHE) approach
\cite{EbrBal:00,FarFer:10,FarFer:11,RaoRaw:03,Rao:00}. Sensor
networks can be employed in many applications, such as monitoring,
exploration, surveillance or tracking  targets over specific
regions. We consider the concept of MHE, as this framework offers
multiple advantages: since a particular minimization problem must be
solved on-line at each step, the observer is optimal with respect to
the associated cost, and moreover, constraints on the state and on
the noise can be taken into account
\cite{FarFer:10,FarFer:11,RaoRaw:03,Rao:00}.

The state estimation problem can be posed  as follows. We assume
that each sensor in the network measures some variables of a
process, computes a local estimate of the entire state of the
system, and exchanges  the computed estimates with its neighbors.
The solution to the estimation problem consists in finding a
methodology which guarantees that all sensors asymptotically reach a
reliable estimate of the overall state of the system. For the
observed process we consider the following nonlinear dynamics:
\[
 x_{t+1}=\phi(x_t) + w_t,
\]
where $x_t \in X  \subseteq \rset^n$ is the state vector and  $w_t
\in W \subseteq \rset^{n}$ represents a white noise with covariance
equal to $Q$.  We also assume  that the sets $X$ and $W$ are convex.
The initial condition  $x_0$ is a random variable with mean $\hat
x_0$ and covariance $\Pi_0$. Measurements on the state vector are
performed by $M$ sensors\footnote{Throughout the paper we will use
the convention that every superscript indicates a sensor/subsystem
index.}, according  to the following sensing model:
\[
y_{t}^i= \theta^i (x_t) +  v_t^i,  \quad  \forall  i=1,\cdots ,M
\]
where $v_t^i \in \rset^{p_i}$ represents white noise  with
covariance matrix $R_i$.   The functions $\phi$ and $\theta^i$ can
be in general nonlinear.

For a given estimation horizon $N \geq 1$, at time $k$ given the
past measurements   $ y_{k-N}^i,\cdots, y_k^i$ provided by the $i$th
sensor and the estimate $\hat x_{k-N}$, we formulate the moving
horizon estimation (MHE) at $k$ as the solution to the following
optimization problem \cite{FarFer:11,RaoRaw:03,Rao:00}:
\begin{align}
 \min_{x_{k-N}, w_{t}}  &  \sum_{i=1}^M  \sum_{t=k-N}^{k}
 || v_t^i ||_{R_i^{-1}}^2 +\!\! \sum_{t=k-N}^{k-1}
|| w_t||_{Q^{-1}}^2  \! +  || x_{k-N} - \hat
x_{k-N}||_{\Pi_{k-N}^{-1}}^2
 \label{est:centralized}  \\
& \text{s.t.}: ~  x_{t+1} = \phi(x_t) + w_t,  \tag{\ref{est:centralized}.1}
\label{est:centralized_dynamics}\\
 &  \qquad  x_t  \in X, \; w_t  \in W   ~~ \forall t,  ~~\tag{\ref{est:centralized}.2}
  \label{est:centralized_set_constr}
\end{align}
where the matrix $\Pi_{k-N}$ is computed recursively from a Riccati
difference equation in a centralized way  \cite{RaoRaw:03}. For the
liner case, the distributed computation of this matrix  can be done
in many ways: e.g. using the steady-state  MHE formulation  (i.e.
computing off-line $\Pi_\infty$, which is the solution of the
corresponding algebraic Riccati equation) or updating $\Pi_{k-N}$
for all the sensors in the same way (using a common covariance
matrix $R$ for all sensors in the Riccati difference equation
update). For the nonlinear case, the update of $\Pi_{k-N}$ in a
distributed fashion is still an open issue.

Note that $  v_t^i =  y_t^i - \theta^i(x_t)$ and using the dynamics
\eqref{est:centralized_dynamics}, we can write $ \sum_{t=k-N}^{k} ||
v_t^i || _{R_i^{-1}}^2$ as a function depending only on $(x_{k-N},
w_{k-N},\cdots, w_{k-1})$. Therefore, by eliminating the states  in
\eqref{est:centralized} using the dynamics
\eqref{est:centralized_dynamics} and introducing the notations:
\begin{align*}
& \qquad \qquad  \qquad  \xb= [x_{k-N}^T \; w_{k-N}^T \cdots w_{k-1}^T]^T, \\
&  f^i (\xb) = \sum_{t=k-N}^{k}  || v_t^i || _{R_i^{-1}}^2 +
\frac{1}{M}\sum_{t=k-N}^{k-1} || w_t||_{Q^{-1}}^2 + \frac{1}{M} ||
x_{k-N} - \hat x_{k-N}||_{\Pi_{k-N}^{-1}}^2,
\end{align*}
 the MHE problem \eqref{est:centralized} can be recast as an optimization problem
 with decoupled cost but a common decision variable $\xb$
 (\textbf{DCx}):
\[ (\textbf{DCx}):
        \left\{
        \begin{array}{l}
             \min_{\xb}  \sum_{i=1}^M f^i(\xb)\\
             \quad \text{s.t.}: ~~ \xb  \in  \textbf{X}, \\
        \end{array}
        \right.
    \]
where the set $\textbf{X} = X \times W^N$.

We assume that the communication network among sensors is described
by a graph $G = (V, E)$, where the nodes in $V= \{1, \cdots, M \}$
represent the sensors and the edge $(i, j) \in E  \subseteq V \times
V$ models  that sensor $j$  sends information to sensor $i$.  Then,
the main challenge is to provide  distributed algorithms for solving
problem \eqref{est:centralized} or equivalently  (\textbf{DCx})
which guarantee that all  the sensors asymptotically  reach a
reliable estimate of the state variables using the information
exchange model given by the graph $G$.

\begin{example}
\label{ex-est} In the particular case where the state and noise
constraints  $x \in X$ and $w \in W$ are described by linear
inequalities (i.e. $X$ and $W$ are polyhedral sets) and  the
dynamics of the process and of the sensors are linear, i.e.
\begin{align*}
& x_{t+1}=A x_t+w_t, \\
 &  y_{t}^i =C_{i}  x_t+v_t^i,   \quad \forall  i=1,\cdots ,M,
\end{align*}
the MHE problem \eqref{est:centralized} can be recast as a
\textit{separable convex quadratic program} with decoupled cost but
a common decision variable in the form (\textbf{DCx}):
\begin{align}
\label {dcx-qp}
& \min_{\xb}   \sum_{i=1}^M  \xb^{T} H_i \xb + q_{i}^T \xb \\
&\qquad \text{s.t.} :  \;    \xb  \in \textbf{X}, \nonumber
\end{align}
where the matrices $H_i$ are positive definite and the constraint
set $\textbf{X}$  becomes in this case polyhedral (described only by
linear inequalities).
\end{example}

%%%%%%%%%%%%%%  DISTRIBUTED CONTROL PROBLEM  %%%%%%%%

\subsection{Distributed optimal control problem}
\label{subsec_optimal} The application that we will discuss  in this
section is the distributed control of large-scale networked systems
with interacting subsystem dynamics, which can be found in a broad
spectrum of applications ranging from traffic networks,  wind farms,
to interconnected chemical plants. Distributed control is promising
in applications for complex systems,  since this framework allows
 us to design local subsystem-based controllers that take into
 account the interactions between different subsystems and physical
constraints.

We consider discrete-time systems which can be decomposed into $M$
subsystems described by difference equations of the form:
\begin{equation}
\label{eq:sub_dnm} x^i_{t+1}  = \phi^i(x^j_t, u^j_t; j  \in
{\mathcal N}^i),   \;\; \forall i=1,\cdots,M,
\end{equation}
where  $x^i_t \in \rset^{n_{i}}$  and   $u^i_t  \in  \rset^{m_{i}}$
represent the state and the input of the $i$th subsystem. The index
set ${\mathcal N}^i$ contains the index $i$ and all the indices of
the subsystems which interact with the subsystem $i$.  We  also
assume that the input and state sequences  must satisfy local
constraints:
\begin{align}
\label{constraints}
 x_{t}^i  \in X^i, \;\;  u_{t}^i \in U^i,  \quad \forall i=1,\cdots,M,  \;\;  \forall t  \geq 0,
 \end{align}
where  the constraint sets $X^i  \subseteq  \rset^{n_i}$ and $U^i
\subseteq \rset^{m_i}$ are usually  compact sets. The system
performance over a prediction horizon of length $N$ is expressed
through a stage cost and a final cost, which are composed of
individual costs for each subsystem $i$ and have the form:
\[    \sum_{t=0}^{N-1} \ell^i(x_t^i,u_t^i) + \ell^i_f(x^i_N).  \]
The centralized optimal control problem over a prediction horizon  $N$ reads:
\begin{align}
 \min_{x_t^i,u_t^i}  & \sum_{i=1}^M \sum_{t=0}^{N-1} \ell^i(x^i_t,u^i_t) + \sum_{i=1}^M \ell_f^i(x_N^i)
 \label{eq:centralized}  \\
& \text{s.t.}: ~  x_0^i = x^i, ~~ x_{t+1}^i = \phi^{i}(x^j_t, u^j_t; j \in {\mathcal N}^i),  \tag{\ref{eq:centralized}.1}
\label{eq:centralized_dynamics}\\
 &  \quad  x^i_t \in X^i, \; u^i_t \in U^i,   ~~ \forall t, i,  ~~\tag{\ref{eq:centralized}.2}
  \label{eq:centralized_set_constr}
\end{align}
where $x^i$ are the values of the initial state for subsystem $i$.
Note that a similar formulation of distributed control for coupled
subsystems with decoupled costs has been given in
\cite{CamJia:02,Dun:07,NegSch:08,NecSav:09,NecSuy:08} in the context
of distributed model predictive control.

Now, we show that the   optimization problem \eqref{eq:centralized}
can be recast as a \textit{separable optimization problem} with a
particular structure. To this purpose, we denote with $\textbf{X}^i
= (X^{i })^{N} \times (U^{i})^{N}$  and
\begin{align*}
& \xb^i=[x_1^{iT} \cdots x_N^{iT} \ u_0^{iT} \cdots u_{N-1}^{iT}]^T, \\
& f^i(\xb^i) = \sum_{t=0}^{N-1} \ell^i(x^i_t,u^i_t)+ \ell_f^i(x_N^i).
\end{align*}
With these notations, problem \eqref{eq:centralized} now reads as an
optimization problem with decoupled cost and sparse coupled
constraints (\textbf{DCCC}):
\[ \textbf{(DCCC) :}
            \left\{
                \begin{array}{l}
\min_{\xb^1,\cdots,\xb^M }  \sum_{i=1}^{M}  f^i(\xb^i) \\
\quad  \text{s.t.}: ~\xb^i \in \textbf{X}^i, ~h^i(\xb^j; ~j \in
\mathcal N^i) =0  ~~\forall i, \\
\end{array}
            \right.
        \]
where the coupled constraints $h^i(\xb^j;~j\in\mathcal N^i) = 0$ are
obtained  from the coupling  between the subsystems, i.e. by
stacking the constraints \eqref{eq:centralized_dynamics} for a given
$i$.

The centralized  optimization problem \eqref{eq:centralized} or
(\textbf{DCCC}) becomes interesting if the computations can be
distributed  among the subsystems (agents), can be  done in parallel
and the amount of information that the agents must exchange is
limited. In comparison with the centralized approach, a distributed
strategy offers a series of advantages: first,  the numerical effort
is considerably smaller since we solve low dimension problems in
parallel and secondly  such a design  is modular, i.e. adding or
removing subsystems does not require any controller redesign.

\begin{example}
\label{dist-ex}
 Many networked systems,  e.g.   wind farms
\cite{MadRan:10}, interconnected chemical processes
\cite{NecSuy:08,Ven:06}, or urban traffic systems \cite{OliCam:10},
can be decomposed into $M$ appropriate linear subsystems:
\begin{align}
\label{sub_dnm}
x^i_{t+1} = A_i x^i_t + B_i u^i_t + \sum_{j \in {\mathcal N}^{-i}}  A_{ij} x^j_t + B_{ij} u^j_t, ~~~  \forall i=1,\cdots, M,
\end{align}
where the index set  ${\mathcal N}^{-i} = {\mathcal N}^{i} - \{i\}
$, i.e. it contains all the indices of the subsystems which interact
with the $i$th subsystem.  If we introduce an auxiliary variable
$w^i_t \in \rset^{p_i}$ to represent the influence of the
neighboring subsystems on the $i$th subsystem  (in applications we
usually have $p_i << n_i$), we can rewrite the dynamics
\eqref{sub_dnm}  as:
\[    x^i_{t+1} = A_i x^i_t + B_i u^i_t + E_i w^i_t, \quad \forall i,      \]
where the matrices  $E_i$ are of appropriate dimensions and
\[ w^i_t   =  \sum_{j  \in  {\mathcal N}^{-i}}  A_{ij}^- x^j_t + B_{ij} ^- u^j_t,   \]
with  the matrices $A_{ij}^-, B_{ij}^-$ being obtained from the
matrices $A_{ij}, B_{ij}$ by removing the rows with all entries
equal to zero.  We consider a quadratic performance index  for each
subsystem $i$ of the form:
\[
\sum_{t=0}^{N-1}  \left (  ||x_t^{i}||^2_{Q_i} + ||u_t^{i}||^2_{R_i}  \right ) + ||x^{i}_N||^2_{P_i},
 \]
where the matrices $Q_i, R_i$ and $P_i$ are positive semidefinite.
  We also assume that the sets $X^i$  and $U^i$ that define the state
   and input constraints \eqref{constraints} are polyhedral.  The centralized control
    problem over the prediction horizon $N$ for this application can be formulated as  follows:
\begin{align}
 \min_{x_t^i, u_t^i, w_t^i}  & \sum_{i=1}^M \sum_{t=0}^{N-1}  ||x_t^{i}||^2_{Q_i} + ||u_t^{i}||^2_{R_i} + ||x^{i}_N||^2_{P_i}
 \label{dc:centralized}  \\
& \text{s.t.}: ~  x_0^i = x^i, ~~ x_{t+1}^i =  A_i x^i_t + B_i u^i_t + E_i w^i_t,  \tag{\ref{dc:centralized}.1}
 \label{dc:centralized_dynamics}\\
& \quad   w^i_t   =  \sum_{j  \in  {\mathcal N}^{-i}}  A_{ij}^- x^j_t + B_{ij} ^- u^j_t, \tag{\ref{dc:centralized}.2}
 \label{dc:centralized_coupling}\\
 & \quad x^i_t \in X^i, \; u^i_t \in U^i   ~~ \forall t,  i. ~~\tag{\ref{dc:centralized}.3}
  \label{dc:centralized_set_constr}
\end{align}
 We can eliminate the state
 variables in the optimization problem  \eqref{dc:centralized}  using the dynamics
  \eqref{dc:centralized_dynamics}. In this  case we can define $\xb^i =[w_0^{iT} \cdots w_{N-1}^{iT} \  u_0^{iT} \cdots u_{N-1}^{iT}]^T$.  Then, the control problem
    \eqref{dc:centralized}  can be recast as a  \textit{separable convex quadratic program} with decoupled cost and coupled constraints in the form (\textbf{DCCC}):
    \begin{align}
\label {dccc-qp}
& \min_{\xb^1,\cdots,\xb^M}  \sum_{i=1}^M  \xb^{iT} H_i \xb^i + q_{i}^T \xb^i \\
&\qquad \text{s.t.} :  \;\;\;   \xb^i  \in \textbf{X}^i, \;\;  \sum_{i=1}^M G_i \xb^i  = g,  \nonumber
\end{align}
where the matrices  $H^i$ are positive semidefinite,  the local
constraint  sets $\textbf{X}^i$ are polyhedral   and the coupled
constraints $\sum_{i=1}^M G_i \xb^i  = g$ are obtained  from the
coupling  between the subsystems, i.e. by stacking
 the constraints \eqref{dc:centralized_coupling} for all $i, t$. Note that the number of rows of the matrices $G^i$ are equal to $N \sum_{i=1}^M p_i$.
\end{example}

%%%%%%%%%%%%%%  COOPERATIVE CONTROL   %%%%%%%%%%%%%%%%%%

\subsection{Cooperative control problem of dynamically uncoupled systems}
\label{subsec_cooperative} Cooperative control for dynamically
uncoupled systems arises in a wide variety of applications like
formation flying, mobile sensor networks, rendezvous problems or
decentralized coordination. The cooperative control problem for
dynamically uncoupled agents consists in controlling a group of
independent subsystems (i.e. with decoupled dynamics), but sharing a
common goal (see e.g. \cite{Dun:07cdc, KevBor:06,KevVer:08}).

We consider a set of $M$ identical subsystems,  having the following
state-space description:
 \[ x^i_{t+1} = \phi(x^i_t,u^i_t),
~~~y^i_t=\theta(x^i_t),~~~\forall i=1,\cdots, M, \]
 where $x^i_t \in
\rset^n$ is the state vector, $u^i_t \in \rset^m$ is the input
vector and $y^i_t \in \rset^p$ is  the  output vector of
 subsystem $i$.  As in the previous section we assume state and input
constraints  of the form \eqref{constraints}.   In the formulation
of cooperative control for uncoupled systems the dynamics of
subsystems are independent from each other, but they share a common
goal. This calls for the minimization of a \textit{cost function}
which involves the states and inputs of each subsystem and their
neighbors as well. In this case we introduce a stage cost at time
$t$ of the form $\ell(x^1_t,\cdots,x^M_t,u^1_t,\cdots,u^M_t)$ and a
final cost $\ell_f(x_N^1,\cdots,x_N^M) $.

The cooperative control problem over a finite horizon of length $N$,
given the initial condition $x^i$ for each subsystem $i$, is
formulated as follows:
\begin{align}
\label{eq_coopcontrol}
 &\min_{x_t^i,u_t^i}  \sum_{t=0}^{N -1}
\ell(x_t^1,\cdots,x_t^M,u_t^1,\cdots,u_t^M) + \ell_f(x_N^1,\cdots,x_N^M) \nonumber \\
& \quad  \text{s.t.}: \; x_0^i=x^i, \;\; x_{t+1}^i = \phi(x^i_t, u^i_t),   \\
& \quad \quad \quad \quad x_t^i \in X^i, \;  u_t^i \in U^i,  \;\; \forall i, t. \nonumber
\end{align}

 Now, let us denote:
\begin{align*}
& \xb^i =[x_1^{iT} \cdots x_N^{iT} \ u_0^{iT} \cdots u_{N-1}^{iT}]^T, \\
& f(\xb^1,\cdots,\xb^M) = \sum_{t=0}^{N-1}  \ell(x_t^1,\cdots,x_t^M,u_t^1,\cdots,u_t^M) + \ell_f(x_N^1,\cdots,x_N^M),
\end{align*}
and $\textbf{X}^i$ the constraint set defined by the state and input
constraints  \eqref{constraints} and by  the $i$th subsystem
dynamics  $x_{t+1}^i = \phi(x^i_t, u^i_t) $ over the prediction
horizon. Using these notations, the previous cooperative control
problem can be recast as an optimization problem with coupled cost
and decoupled constraints (\textbf{CCDC}):
\[ \textbf{(CCDC) :}
            \left\{
                \begin{array}{l}
                    \min_{\textbf{x}^1,\cdots,\textbf{x}^M} f(\textbf{x}^1,\cdots,\textbf{x}^M)\\
                    \qquad \text{s.t.}: \; \textbf{x}^i \in \textbf{X}^i.   \\
                \end{array}
            \right.
        \]
We are interested in finding efficient parallel algorithms for solving problem (\textbf{CCDC}).

\begin{example}
\label{ex_coop} We consider the formation flying  for a group of
satellites that are distributed  along a circular orbit with
independent dynamics but they have to maintain a constant distance
with respect to the two nearest neighbors (see e.g.
\cite{KevVer:08}).  Using a discretized version of the linear
Clohessy-Wiltshire equations  of the $i$th satellite for a nominal
circular trajectory  \cite{Kap:76}:
\begin{equation*}
\left\{
\begin{array}{l}
      \ddot{x}^{1,i}=3\omega_n^2x^{1,i}+2\omega_n\dot{x}^{2,i}+a^{1,i}\\
     \ddot{x}^{2,i}=-2\omega_n\dot{x}^{1,i}+a^{2,i}\\
     \ddot{x}^{3,i}=-\omega_n^2x^{3,i}+a^{3,i},\\
\end{array}
  \right.
  \end{equation*}
where $x^{1,i}$, $x^{2,i}$, $x^{3,i}$ are  the displacements in the
radial, tangential and out-of-plane direction, $a^{1,i}$, $a^{2,i}$,
$a^{3,i}$ represent the accelerations of the satellite $i$ due to
propulsion or external disturbances and $\omega_n$ is the angular
velocity at which the orbit is covered, we obtain a discrete-time
linear system for the $i$th satellite of the form
\begin{equation*}
\left\{
\begin{array}{l}
     x_{t+1}^i=Ax_{t}^i + Bu_{t}^i\\
     y_{t}^i=Cx_{t}^i, \\
\end{array}
  \right.
  \end{equation*}
with    $x^i_t  \in \rset^6$ and $u^i_t=[a^{1,i}_t \ a^{2,i}_t \
a^{3,i}_t]^T \in \rset^3$ being  the state, respectively the input
vectors of satellite $i$ and  we consider as output
$y^i_t=[x^{1,i}_t \  x^{2,i}_t  \  x^{3,i}_t]^T$,  the vector of
absolute positions of the satellite.  We also assume input
constraints of the form:
\[ u_{\min}  \leq u_t^i \leq u_{\max} \quad \forall i, t.   \]
Since the goal is to maintain a constant distance with respect to
the two nearest neighbors, we choose the following stage cost at
time $t$:
    \begin{align*}
    & \ell(x^1_t,\cdots,x^M_t,u^1_t,\cdots,u^M_t)=  \sum_{i=1}^M  ||2y^i_t - y^{i+1}_t - y^{i-1}_t||_{Q_i}^2+||u^i_t||_{R_i}^2,
    \end{align*}
where $Q_i, R_i$ are positive definite matrices. We assume the final
cost $\ell_f = 0$.  Despite the fact that the output $y^i$
represents the absolute positions of the $i$th satellite, using the
stage cost from above, the formation flying becomes a problem based
on relative positions between the satellites instead of the absolute
ones. In this case the cooperative control problem
\eqref{eq_coopcontrol} over a finite horizon $N$  can be recast as a
\textit{convex quadratic problem} with coupled cost and decoupled
constraints in the form (\textbf{CCDC}) :
%\begin{equation} \left\{
\begin{align}
\label {ccdc-qp}
& \min_{\xb^1,\cdots,\xb^M}
 \left[ \begin{array}{c}
\xb^1\\
\vdots \\
\xb^M
\end{array} \right] ^T
%\left[ \begin{array}{ccccccccccc}
%* & * & * & 0 & 0 & \cdots & 0 &  0 &  0 & * &  * &\\
%* & * & * & * & 0 & \cdots & 0 &  0 &  0 & 0 &  * &\\
%* & * & * & *&  *& \cdots & 0 &  0 &  0 & 0&  0 &\\
%\ddots & \ddots & \ddots & \ddots & \ddots & \ddots & \ddots &  \ddots &  \ddots & \ddots &  \ddots &\\
%0 & 0 & 0 & 0 & 0 & \cdots & * &  * &  * & * &  * &\\
%* & 0 & 0 & 0 & 0 & \cdots & 0 &  * &  * & * &  * &\\
%*& * & 0 & 0 & 0 & \cdots & 0 &  0 &  * & * &  * &
%\end{array} \right]
 \left[ \begin{array}{c}
H_{ij}
\end{array} \right] _{ij}
\left[ \begin{array}{c}
\xb^1\\
\vdots \\
\xb^M
\end{array} \right]  +
\left[ \begin{array}{c}
q_1\\
\vdots \\
q_M
\end{array} \right] ^T
 \left[ \begin{array}{c}
\xb^1\\
\vdots \\
\xb^M
\end{array} \right]  \\
& \qquad \text{s.t.} : \; \xb^i  \in \textbf{X}^i, \nonumber
\end{align}
%\right \end{equation}
where the blocks of the positive semidefinite  Hessian matrix $H =
[H_{ij}]_{ij}$  satisfies $H_{ij} = 0$ if $|i-j| > 3$ for all $i, j$
and the sets $\textbf{X}^i$ are polyhedral.

\begin{remark}
\label{rem-elimx}
(i)   Note  that we can eliminate the states $x_1^i,\cdots,x_N^i$ using the
 dynamics of the $i$th satellite and keeping  only the inputs over the prediction
 horizon as decision variables, i.e. we
  may redefine  $\xb^i =[\ u_0^{iT} \cdots u_{N-1}^{iT}]^T$.
   In this case $H$ becomes positive definite and the sets  $\textbf{X}^i$
    are described only by linear inequalities. \\
(ii)  In many  applications we can move the coupling terms from the
cost to the constraints  by introducing auxiliary variables, i.e we
can recast an optimization problem with coupled cost but decoupled
constraints (\textbf{CCDC}) to one with decoupled cost but coupled
constraints (\textbf{DCCC}).   E.g., in our satellite formation
application we can define the coupling constraints $w_t^i =
y_t^{i-1} + y_t^{i+1} $ and then we can associate a local stage cost
for each satellite $i$ as $\ell^i(x_t^i,w_t^i,u_t^i) = ||2 C x^i_t -
w^{i}_t ||_{Q_i}^2+||u^i_t||_{R_i}^2$ but with coupled dynamics
$w_t^i = C (x_t^{i-1} + x_t^{i+1})$.   We can also do the other way
around: we can reformulate a  (\textbf{DCCC}) into a (\textbf{CCDC})
problem (e.g. by moving the coupling constraints
\eqref{eq:centralized_dynamics} into the cost, see Section
\ref{subsec_cooperative_coupled}). Depending on applications one
formulation might be preferred against the other (see also Section
\ref{subsec_cooperative_coupled} below).
\end{remark}

\end{example}

%%%%%%%%%%%%%%%%%%%%%%%%%%%%%%%%%%%%%%%%%%%%%%%%%%%%%%%%%%%%%

\subsection{Cooperative control problem of dynamically coupled systems}
\label{subsec_cooperative_coupled}

In this section we discuss  the cooperation-based optimal control
problem for a group of dynamically coupled  subsystems
\cite{CamOli:09,MaePen:10,PanWri:09,RawSte:08,VenRaw:05,Ven:06}. For
the $i$th subsystem we consider the following linear dynamics:
\begin{align}
\label{sub_dnm2}
  &x^i_{t+1} = A_i x^i_t + B_i u^i_t + \sum_{j \in {\mathcal N}^{-i}} B_{ij} u^j_t,~~~\forall i=1,\cdots, M.
\end{align}
Note that the dynamics described in \eqref{sub_dnm2} are a
particular case of \eqref{sub_dnm}. We also assume local input
constraints $u^i_t \in U^i$, where $U^i$ are convex sets.

For each subsystem we define a local stage cost $\ell^i(x^i,u^i)$
and a  terminal cost $\ell^i_f(x^i)$. The local cost for each
subsystem on a finite horizon of length $N$ will be of the following
form:
\begin{equation}
f^i(\overline{\xb}^i, \overline{\ub}^i) = \sum_{i=0}^{N-1}
\ell^i(x^i_t, u^i_t) + \ell^i_f(x^i_N),
\end{equation}
where we denote with
\begin{equation} \overline{\xb}^i=[x^{iT}_1\cdots
x^{iT}_N]^T,~~~\overline{\ub}^i=[u^{iT}_0\cdots u^{iT}_{N-1}]^T.
\end{equation}

In order  to provide a cooperative behavior between subsystems we
replace each local cost $f^i$ with one that represents the
systemwide impact of local control actions. One choice is to employ
a strong convex combination of local subsystems' costs as the global
objective function for the entire system. In these conditions, the
cooperative control problem  for coupled systems on a finite horizon
$N$ will have the form:
 \begin{align}
 \min_{\overline{\xb}^i,\overline{\ub}^i}  & \sum_{i=1}^M \alpha_i f^i(\overline{\xb}^i,\overline{\ub}^i)
 \label{dc:coop}  \\
& \text{s.t.}: ~~ x_{t+1}^i =  A_i x^i_t + B_i u^i_t +
\sum_{j \in {\mathcal N}^{-i}} B_{ij} u^j_t, \;\;\;\;  x_0^i = x^i,  \tag{\ref{dc:coop}.1}  \label{dc:coop_dynamics}\\
 &  ~~u^i_t \in U^i   ~~ \forall t,  i, ~~\tag{\ref{dc:coop}.2}
  \label{dc:coop_set_constr}
\end{align}
where $\alpha_i>0$  and sum to 1. Note that in this form problem
\eqref{dc:coop} is a particular case of problem (\textbf{DCCC}),
where the variables associated to the $i$th subsystem are given by
$[\overline{\xb}^{iT} \;  \overline{\ub}^{iT}]^T$. However, by
eliminating the states in  \eqref{dc:coop} using the global dynamic
model  (\ref{dc:coop}.1) we obtain a coupled objective function in
the local variables $\xb^i = \overline{\ub}^i$ (i.e. in the local
control actions) and decoupled constraints,  which is  a particular
case of (\textbf{CCDC}) problem (see also Remark
\ref{rem-elimx}(ii)).

%%%%%%%%%%%%%% DECOMPOSITIONS METHODS  %%%%%%%%%%

\section{Parallel and distributed optimization algorithms  for solving coupled optimization problems}
\label{sec_decomposition}

In this section we present several parallel and distributed
algorithms for solving the optimization problems  arising in
applications from estimation and control discussed  in Section
\ref{sec_application} and analyze their properties and performances,
in particular we define conditions  for which these algorithms
converge\footnote{For simplicity of the exposition, in this section
we assume that all the functions are differentiable.}. The presented
algorithms  can be classified, on the one hand in ``centralized''
algorithms (that in general take advantage of the sparsity of the
problem and solve in parallel low dimension subproblems)  and
distributed algorithms (that take into account explicitly
information restrictions in the network  and combine consensus
negotiations with optimization methods to solve distributively the
problem) and on the other hand in primal and dual decomposition
algorithms.  The first class is based on  decomposing the original
optimization problem, while the second consists in decomposing the
corresponding dual  problem.

For a given problem representation there are often many choices of
distributed algorithms, each with possible different
characteristics: e.g.  rate of convergence, tradeoff between local
computation and global communication, and quantity of message
passing. Which alternative is the best depends on the specifications
of the application. However, for each algorithm  we will discuss in
details their main characteristics in terms of performance and
properties.

\subsection{Distributed  gradient algorithms for  optimization problems of type (\textbf{DCx})}
\label{sec-dgp} In this section  we study  several distributed
algorithms  for solving  separable optimization problems with
decoupled cost and common decision variables in the form
(\textbf{DCx}), that e.g.  appear in the context of state estimation
in  sensor networks (see Section \ref{subsec_estimation}).   We
associate to the set  of agents (e.g. sensors)  a graph $G=(V, E)$
and then such distributed algorithms must satisfy  the following
constraint:  the computations will be performed on all nodes in
parallel, and the communication between nodes is restricted to the
edges of the graph.  Distributed optimization algorithms are mainly
based on combining consensus negotiations (as an efficient method
for information fusion) with optimization methods  \cite{JohKev:08,
NecDum:10, NedOzd:09, NedOzd:10,ZhuMar:10} to solve  distributively
problems of type (\textbf{DCx}).

First we introduce  the consensus problem for a group of  $M$ agents
that considers conditions under which using a certain
message-passing protocol, the local variables of each agent will
converge to the same value  \cite{Mor:09,OlfFax:07,XiaWan:08}. There
exist several results related to the convergence of local variables
to a common value using various information exchange protocols among
agents \cite{OlfFax:07,OlsTsi:06,XiaWan:08}.   One of the most used
models for consensus is based on the following discrete-time
iteration: to generate an estimate at iteration $k+1$, agent $i$
forms a convex combination  of its estimate $\xb^i_k$ with the
estimates received from other agents:
\begin{align*}
& \xb^i_{k+1} = \sum_{j=1}^M  \gamma^{ij}_k  \xb^j_k,
\end{align*}
where $\gamma^{ij}_k$ represent nonnegative
weights\footnote{Naturally, an agent $i$  assigns zero weight  to
the estimates $\xb^j$ for those agents $j$ whose estimate
information is not available at the update time.}  satisfying $\sum_
j \gamma_k^{ij} = 1$.  At each iteration $k$ the information
exchange among agents  can be represented by a graph  $(V, E_k)$,
where $E_k =\{ (i, j): \gamma_{k}^{ij} >0 \}$.  We can also
introduce the graph $(V, E_\infty)$, where $E_\infty = \{ (i, j):
(i, j) \in E_k \; \text{for infinitely many} \; k \}$. The graphs
$(V, E_k)$ satisfy the \textit{bounded  interconnection interval
property} if there exists an integer $\tau$  such that for any $(i,
j) \in E_\infty$ agent $j$ sends its information to agent $i$ at
least once every $\tau$  consecutive iterations. It has been
 proved in \cite{NedOzd:09} that under certain assumptions on the weights $\gamma^{ij}_k$
 (e.g.  stochasticity  of the matrix $\Gamma_k=[\gamma^{ij}_k]_{ij}$,
 strong connectivity property of $(V, E_\infty)$ and
 bounded interconnection interval property),
 the states $\xb^i_k$ of all agents converge to the same state
 $x^*$. Similar convergence results can be found in
 \cite{Mor:09,XiaWan:08}.

We return now to our optimization problem  of type (\textbf{DCx}).
In \cite{NedOzd:10}  a distributed projected gradient algorithm is
analyzed, which basically combines  the consensus iteration
presented above with a projected gradient update to generate the
next estimate of the optimum.   More specifically, an agent $i$
updates its estimate by combining the estimates received from its
neighbors, then  taking a gradient step to minimize its objective
function $f^i$ and finally projecting on the set $\textbf{X}$:
\begin{algorithm}[h!]
\textbf{Algorithm  dgp1}
 \[   \vb^i_k =\sum_{j=1}^M \gamma^{ij}_k \xb^j_k, \quad
       \xb^i_{k+1} =  \left [  \vb^i_k  - \alpha_k  \nabla  f^i(\vb^i_k)  \right ]_{\textbf{X}}   \]
\end{algorithm}

where $\alpha_k$ is a  common step size, $ \nabla  f^i$ denotes the
gradient of the function $f^i$,  and
 $[ \cdot ]_{\textbf{X}}$ denotes the Euclidian projection on the set $\textbf{X}$.
 The following convergence result holds for Algorithm \textbf{dgp1} :
\begin{theorem}  \cite{NedOzd:10}
\label{th-dgp} For the optimization problem  (\textbf{DCx})  we
assume that all the functions $f^i$ are convex and have bounded
gradients, the set $\textbf{X}$ is convex and the step size
satisfies $\sum_k \alpha_k = \infty$ and $\sum_k \alpha_k^2 <
\infty$. Moreover, we assume that the weights $\gamma_k^{ij}$
satisfy the following properties: the matrices $\Gamma_k =
[\gamma_k^{ij}]_{ij}$ are doubly stochastic, the graph $(V,
E_\infty)$ is connected and the bounded interconnection interval
property holds.  Then, the distributed projected gradient Algorithm
\textbf{dgp1} converges to an optimum of problem (\textbf{DCx}).
\end{theorem}

An  interesting variant  of  a distributed gradient   projected
algorithm has been provided in \cite{JohKev:08}. Compared to the
previous distributed gradient Algorithm \textbf{dgp1}, in
\cite{JohKev:08} a fixed connected graph $(V, E)$ is taken over all
iterations and the information exchange among the agents is
represented by a doubly stochastic matrix $\Gamma =
[\gamma^{ij}]_{ij}$ such that $\gamma_{ij} > 0$ if $(i, j) \in E$.
In this algorithm, first each agent implements the gradient update
locally and then it  runs a number $\mu$ of consensus iterations
with its neighbors:
\begin{algorithm}
\textbf{Algorithm  dgp2}
 \[
       \xb^i_{k+1} = \left [ \sum_{j=1}^M   \Gamma^{\mu}_{ij}
         \left ( \xb^j_k   - \alpha_k  \nabla  f^j(\xb^j_k) \right )
         \right ]_{\textbf{X}}   \]
\end{algorithm}\\
where $\Gamma^{\mu}_{ij} $  denotes the $(i, j)$ entry of the matrix
$\Gamma^{\mu}$.  Under similar assumptions as in Theorem
\ref{th-dgp}, the authors in \cite{JohKev:08} proved convergence of
Algorithm  \textbf{dgp2} for a constant step size and for a
sufficiently large $\mu$.

In the case when the set $\textbf{X}$ is explicitly defined through
a finite set of equalities and inequalities, an algorithm based on a
penalty primal-dual approach has been recently proposed in
\cite{ZhuMar:10}. This algorithm allows the agents exchange
information over networks with time-varying topologies and
asymptotically agree on an optimal solution and the optimal value.

Another interesting approach  for solving the optimization problem
(\textbf{DCx}), but in a serial fashion,  can be found in
\cite{NedBer:01} where an incremental gradient method is presented.
Each step of the algorithm is a gradient iteration for a single
component function $f^i$, and there is one step per component
function. Thus, an iteration can be viewed as a cycle of $M$
subiterations, so that at $k+1$:
\begin{align*}
& \xb_{k+1} = z_{M,k}, \quad z_{0,k}=\xb_k, \\
& z_{i,k} = \left [z_{i-1,k} - \alpha_k \nabla  f^i(z_{i-1,k})
  \right ]_{\textbf{X}} ~~~~ \forall i=1,\cdots,M.
\end{align*}
For convex problems, using an appropriate step size $\alpha_k$, the
authors in \cite{NedBer:01} show that this algorithm has much better
practical rate of convergence  than the classical gradient method.

\begin{remark}
(i)  The convexity assumptions on the functions $f^i$ and the set
$\textbf{X}$  for convergence of the two Algorithms \textbf{dgp1}
and \textbf{dgp2} are usually satisfied in many applications: see
e.g. the state estimation problem for linear systems
discussed in Example \ref{ex-est} which leads to the convex quadratic program \eqref{dcx-qp}. \\
(ii)  One of the main challenges when solving problems of type
(\textbf{DCx}) is the time-dependent communication topology, as
communication links can change due to changing distances, obstacles,
or disturbances. While in  \cite{JohKev:08} a constant topology is
assumed for Algorithm \textbf{dgp2}, the Algorithm \textbf{dgp1} and
the algorithm from \cite{ZhuMar:10} are based on a changing
topology, which makes them more suitable in practical applications.
Moreover, the cyclical incremental algorithm \cite{NedBer:01} can be
implemented  only when each agent
 identifies a suitable downstream and upstream neighbor. Note the existence
 of a cycle is a stronger  assumption than connectivity. \\
(iii) From simulations we have observed that the algorithms from
\cite{JohKev:08,NedOzd:10,ZhuMar:10} are very sensitive to the
choice of the weights that must be tuned, since they are considered
as parameters in these methods. These algorithms do not provide a
mathematical way of choosing the weights from the consensus
protocol, which has a very strong influence on the convergence rate
of these methods. Recently in \cite{NecDum:10}, a distributed
algorithm  has been derived for solving particular cases of problems
of type (\textbf{DCx}),  where the nonnegative weights corresponding
to the consensus process are interpreted as dual variables and thus
they are updated using arguments from duality theory. Moreover, if
the network is not densely connected (i.e. each sensor has  a large
number of neighbors), one can expect  the performance of these
algorithms from \cite{JohKev:08,NedOzd:10,ZhuMar:10} to be worse
than that of the cyclic incremental gradient \cite{NedBer:01}.
\end{remark}

\begin{table}[h!]
\label{table_est}
\begin{center}
\begin{tabular}{|c|c|c|c|c|}
\hline $M$    & $N$    &       \text{nr. it. dgp1}         &   \text{nr. it. dgp2}  \\
\hline $10$  & $10$   &      $5.627$                                                 &     $586$     \\
\hline $10$  & $20$   &      $8.447$                                                 &     $746$     \\
\hline $20$  & $10$   &      $10.651$                                               &     $1.854$   \\
\hline $20$  & $20$   &      $14.758$                                               &     $2.571$   \\
\hline
\end{tabular}
\end{center}
\caption{State estimation problem  Example \ref{ex-est}:  we
consider $M=10, 20$ sensors, a
 linear system with $5$ states and a prediction horizon $N=10, 20$.   We solve the convex
  quadratic program \eqref{dcx-qp} with the
 accuracy of the solution $\epsilon=10^{-2}$.  We assume fixed weights in both algorithms
 such that $\gamma^{ij}=0$  for $|i-j|>1$  and $\mu=10$.
 From simulations we observe that Algorithm \textbf{dgp2} works better than Algorithm
 \textbf{dgp1} in terms of the number of gradient iterations. However, Algorithm
   \textbf{dgp2} needs to perform for each gradient iteration also $\mu=10$ consensus steps.}
 \vspace{-0.1cm}
\end{table}

%%%%%%%%%%%%%%%%%%%%%%%%%%%%%%%%%%%%%%%%%%%%

\subsection{Decomposition algorithms for solving  optimization problems  (\textbf{DCCC})}

In this section  we present several decomposition  algorithms for
solving  separable optimization problems with decoupled cost but
coupled constraints in the form (\textbf{DCCC}).  Distributed
control  for complex processes  with interacting subsystem dynamics
usually   leads to such optimization problems  (see e.g. Section
\ref{subsec_optimal}).  We discuss two classes of decomposition
principles: primal and dual. We use the terms primal and dual  in
their mathematical programming meaning: primal indicates that the
optimization problems are solved  using the original formulation and
variables  and dual indicates that the original problem has been
rewritten using Lagrangian relaxation.

Compared to the general  formulation of problem  (\textbf{DCCC}), we
  focus in this section on decomposition methods for the particular
case of separable \textit{convex} problems with decoupled cost and
coupled constraints\footnote{For  the nonconvex   case  of  problem
(\textbf{DCCC}) we can still  obtain decomposition algorithms by
combining  sequential quadratic programming or sequential convex
programming, in order to linearize the nonlinear coupled
constraints,  with  decomposition methods that address the
decomposable convex problems  (see e.g.  \cite{NecSav:09}).}:
\[ \textbf{(conv-DCCC):}
            \left\{
                \begin{array}{l}
\min_{\xb^1,\cdots,\xb^M }  \sum_{i=1}^{M}  f^i(\xb^i) \\
\quad  \text{s.t.}: ~\xb^i \in \textbf{X}^i,  ~  \sum_{i=1}^M G_i \xb^i  = g,\\
\end{array}
            \right.
        \]
where  we consider that for all $i$ the coupled constraints
$h^i(\xb^j; ~j \in \mathcal N^i) =0$ in  problem (\textbf{DCCC})
become linear and can be written compactly as $\sum_{i=1}^M G_i
\xb^i  = g $, with $G_i \in \rset^{n_\lambda \times n_{\xb^i}}$. For
simplicity of the exposition  the following  assumptions hold for
problem \textbf{(conv-DCCC)} (for general case of convex problems
see \cite{NecSuy:08, NecSuy:09}):
\begin{assumption}
\label{ass-dual} Each function $f^i$  is convex quadratic  and
$\textbf{X}^i$ are compact convex sets. Moreover,   the Slater's
condition holds, i.e. there exist $\xb^i  \in
\text{int}(\textbf{X}^i)$ such that
 $\sum_{i=1}^M G_i \xb^i  = g$.
\end{assumption}
From Example \ref{dist-ex} we have seen  that centralized optimal
control for interconnected linear systems leads to such a separable
convex quadratic formulation, e.g. \eqref{dccc-qp}.

We begin with primal  decomposition (see e.g.
\cite{BoyXia:03,ChiLow:06,PalChi:06,Sil:72} and the references
therein). We can decompose the original problem (\textbf{conv-DCCC})
as follows:   we introduce some auxiliary variables in order to
separate  the coupled linear equality constraints, i.e. we introduce
the new variables $\tb^1, \cdots, \tb^{M-1}$, and obtain $M$
subproblems:
\begin{align*}
                    (\Pp^i):  \quad  \psi^i (\tb^i) = \min_{\xb^i}  \{ f^i(\xb^i) :  \; \xb^i    \in \textbf{X}^i,   \;  G_i  \xb^i  = \tb^i \}
\end{align*}
for $i= 1, \cdots, M -1$  and  the $M$th     subproblem
    \begin{align*}
&  (\Pp^M):  \;\;  \psi^M  (\tb^1,\cdots,\tb^{M-1}) =  \min_{\xb^M} \{ f^M(\xb^M): \;  \xb^M   \in   \textbf{X}^M, \ \sum_{i=1}^{M-1} \tb^i  \!+\!  G_M \xb^M  \!=\! g \}.
\end{align*}
The separable convex problem (\textbf{conv-DCCC}) reduces to solving
the unconstrained convex  \textit{primal problem}  $(\Pp \Pp)$
\cite{Sil:72}:
\[
(\Pp \Pp): \qquad \min_{\tb^1,\cdots,\tb^{M-1}}   \psi (\tb^1,  \cdots, \tb^{M-1}),
  \]
 where $\psi (\tb^1,  \cdots, \tb^{M-1})  = \psi^1(\tb^1) + \cdots + \psi^{M-1} (\tb^{M-1})
 + \psi^M  (\tb^1,  \cdots, \tb^{M-1})$.  Conditions for well-posedness  of the primal problem  $(\Pp \Pp)$  can be found in  \cite{Sil:72}.
   Let $\xb^i(\tb^i)$ and $\lambda^i(\tb^i)$ be
 the optimal solution and the corresponding optimal Lagrange multiplier for the equality
  constraints $G_i  \xb^i = \tb^i$, respectively,  for subproblem $\Pp^i$ given $\tb^i$,
  with $i= 1,\cdots,M -1$.    Similarly, we define $\xb^M(\tb^1,\cdots,\tb^{M-1})$ and
  $\lambda^M(\tb^1,\cdots, \tb^{M-1})$  for subproblem $\Pp^M$.
  Although the function $\psi$ is potentially non smooth, assuming that Slater's
  condition for the convex problem (\textbf{conv-DCCC}) holds (according to
   Assumption \ref{ass-dual}), the following vector is a
   subgradient\footnote{A vector $s \in \rset^n$ is a subgradient of
    $f: \rset^n \rightarrow  \rset$ at a point $x \in  \text{dom} f$
    if for all $y \in \text{dom} f$ we have  $f(y)\geq  f(x) + s^T (y-x)$.}
     of $\psi$   at $(\tb^1,\cdots, \tb^{M-1})$  \cite{BerNed:03,Sil:72}:
\[  \left [ \lambda^M(\tb^1,\cdots, \tb^{M-1}) - \lambda^1(\tb^1)  \cdots
  \lambda^M(\tb^1,\cdots, \tb^{M-1})  - \lambda^{M-1}(\tb^{M-1})  \right ]^T. \]

\begin{algorithm}[h!]
\textbf{Algorithm primal subgradient (PS)}
\begin{align*}
&\xb^{i}_{k} =  \xb^i (\tb^i_k),  \quad  \lambda^{i}_{k} =  \lambda^i (\tb^i_k) \;\; \text{for} \; i=1,\cdots, M-1 \\
&  \xb^{M}_{k} = \xb^M(\tb^1_k,\cdots,\tb^{M-1}_k),   \quad \lambda^{M}_{k} = \lambda^M(\tb^1_k,\cdots,\tb^{M-1}_k),  \\
&  \qquad \qquad \tb^i_{k+1} = \tb^i_k  - \alpha_k  ( \lambda^M_{k}
- \lambda^i_{k} ),
\end{align*}
\end{algorithm}
where  $\alpha_k$ is a step size.
\begin{remark}
\label{rem_step} The step size $\alpha_k$ can be chosen in two ways:
(i) it can vary but satisfying $\sum_k \alpha_k = \infty$ and
$\sum_k \alpha_k^2 < \infty$; (ii) $\alpha_k$ is constant for all
$k$.
\end{remark}
Under Assumption  \ref{ass-dual} the convergence of this primal
subgradient algorithm is obvious, due to the equivalence between the
(\textbf{conv-DCCC}) problem and the convex primal problem  $(\Pp
\Pp)$. When the primal problem  $(\Pp \Pp)$ (called also the
\textit{master problem}) is solved using this scheme, the method has
an interesting economic interpretation: at each iteration the master
program allocates the resources (by choosing
 $\tb_k^i$) and the nodes return the prices  associated with this choice $\lambda^{i}_{k}$. The iteration continues until the prices have reached the equilibrium.

We now discuss dual decomposition
\cite{BerTsi:89,ChiLow:06,KonLeo:96,NecSuy:09,Spi:85,Tse:91}. In
dual decomposition methods we have the following economic
interpretation: the master problem sets the prices for the resources
to each subproblem which has to decide the amount of resources to be
used depending on the price. The iteration  continues  until the
best pricing strategy is obtained. Clearly, if the coupled
constraints $\sum_i G_i \xb^i = g$ are absent, then the problem
(\textbf{conv-DCCC}) can be decoupled. Therefore it makes sense to
relax these coupled constraints  using duality theory.  We construct
the \textit{partial augmented Lagrangian}:
\begin{align}
    \label{lagrangian}
    L_\mu(\xb, \lambda)  =\sum_{i=1}^M f^i(\xb^i)  + \mu P_{\textbf{X}^i} (\xb^i) +  \lambda^T  (\sum_{i=1}^M  G_i \xb^i  - g),
\end{align}
where $\mu>0$ and the functions   $P_{\textbf{X}^i}$ associated  to
the sets $\textbf{X}^i$  (usually called \textit{prox functions})
must have certain  properties  explained below. We also define the
corresponding \textit{augmented dual function}:
\begin{align}
\label{dip}
d_\mu(\lambda) = \min_{\xb^i  \in \textbf{X}^i}  L_\mu(\xb,\lambda),
\end{align}
and from the structure of $L_\mu$ we obtain that    \eqref{dip}
decouples in $M$ subproblems
\[   \xb^i( \mu,\ld) =  \arg \min_{\xb^i  \in \textbf{X}^i}    f^i(\xb^i)  +  \mu P_{\textbf{X}^i} (\xb^i)  + \lambda^T G_i \xb^i.  \]
 We are interested in the properties of  the family of augmented dual functions
$\{ d_\mu \}_{\mu > 0}$. Note that $\lim_{\mu \to 0} d_\mu (\lambda)
= d_0 (\lambda)$, where $d_0 (\lambda) = \min_{\xb^i \in
\textbf{X}^i} L_0 (\xb, \lambda)$ is the standard dual function,
whenever the prox functions $P_{\textbf{X}^i}$ are chosen to be
continuous on the compact sets $\textbf{X}^i$ or are barrier
functions associated to these sets (see \cite{Nes:04}). The goal is
to maximize the augmented dual function for $\mu$ sufficiently
small:
\[
 \max_\lambda d_\mu (\lambda),
 \]
in order to find an approximation of the optimal Lagrange multiplier
$\lambda^* = \arg \max_\lambda d_0 (\lambda)$    and then to recover
an approximation of the corresponding optimal primal variables
$\xb^{i*}$. We distinguish three algorithms, depending on the choice
of the constant $\mu$ and of the prox functions
 $P_{\textbf{X}^i}$:
\begin{itemize}
\item[(I)]  \textbf{dual subgradient algorithm}:   $\mu = 0$ and $P_{\textbf{X}^i} = 0$
\item[(II)]  \textbf{dual fast gradient algorithm}:   $\mu > 0$ and $P_{\textbf{X}^i} $ are strongly convex functions
\item[(III)]  \textbf{dual interior-point algorithm}:  $\mu > 0$ and $P_{\textbf{X}^i} $ are barrier  functions for the sets
$\textbf{X}^i$.
\end{itemize}
The next theorem provides the main properties of the augmented dual
function:
\begin{theorem}
\cite{NecSuy:08, NecSuy:09}   Under Assumption  \ref{ass-dual}, the augmented dual function $d_\mu$ is characterized as follows:\\
(I)  For any $\mu \geq 0$ and convex functions $P_{\textbf{X}^i} $ a
subgradient of   $d_\mu$ at $\lambda$ is given by  $ \sum_{i} G_i
\xb^i( \mu,\ld) - g$.  (II)  For  $\mu > 0$ and strong convex
functions $P_{\textbf{X}^i} $ the function $d_\mu$  has a Lipschitz
continuous gradient.    (III)  For   $\mu > 0$ and barrier functions
$P_{\textbf{X}^i} $   the function $d_\mu$ is self-concordant.
\end{theorem}
 We denote
$\xb_k^i = \xb^i (\mu_k, \lambda_k)$. The  iterations of the three
algorithms are:
\begin{algorithm}[h]
\textbf{Algorithm  dual subgradient  (DS)}
 \[
\lambda_{k+1}= \lambda_k + \alpha_k (\sum_{i=1}^M G_i \xb^i_k - g)
\]
\textbf{Algorithm dual fast gradient  (DFG)}
\[
\bar \lambda_{k+1}= \lambda_k + \frac{1}{L_{\mu_k}}(\sum_{i=1}^M G_i \xb^i_k - g), \quad \lambda_{k+1}= \bar \lambda_{k+1} + \beta_k (\bar \lambda_{k+1} - \lambda_k)
\]
\textbf{Algorithm dual interior-point (DIP)}
\[
\lambda_{k+1}= \lambda_k + \alpha_k  \left ( \nabla^2 d_{\mu_p}(\lambda_k)  \right )^{-1} \nabla d_{\mu_p} (\lambda_k) \;\; \text{as} \;\; \mu_p  \to 0,
\]
\end{algorithm}\\
where   $\alpha_k$ is a  step-size that can be chosen as in Remark
\ref{rem_step} for algorithm (\textbf{DS}) or satisfying Armijo rule
\cite{Nes:04} for algorithm (\textbf{DIP}), $L_\mu$ is the Lipschitz
constant of the gradient $\nabla d_\mu$ and $\beta_k>0$ is defined
iteratively as in \cite{Nes:04}. Moreover, in the dual
interior-point algorithm (\textbf{DIP}) we have an outer iteration
in $p$ where we decrease $\mu_p  \to 0$ and an inner iteration in
$k$ where we need to generate vectors close to the central path
using Newton updates with $\nabla^2 d_{\mu}(\lambda)$ representing
the Hessian of the augmented dual function $d_\mu$ at $\lambda$ (see
\cite{NecSuy:09} for more details).

The convergence of these three algorithms (\textbf{DS}),
(\textbf{DFG}) and  (\textbf{DIP})   can be established under
suitable assumptions on problem (\textbf{conv-DCCC})  and on  the
prox functions $P_{\textbf{X}^i} $:
\begin{theorem}\cite{NecSuy:08, NecSuy:09}
If Assumption  \ref{ass-dual} holds for the separable convex problem
(\textbf{conv-DCCC}), then all three algorithms (\textbf{DS}),
(\textbf{DFG}) and (\textbf{DIP}) are convergent under a suitable
choice of the step-size. Moreover,  the dual fast gradient algorithm
(\textbf{DFG})  has complexity ${\mathcal O}
(\frac{c_1}{\epsilon})$, while the dual interior-point algorithm
(\textbf{DIP})  has complexity ${\mathcal O} \left ( c_2 \log
(\frac{c_3}{\epsilon}) \right)$, where $\epsilon$ is the accuracy of
the approximation of the optimum for problem (\textbf{conv-DCCC})
and $c_i$ are some positive constants.
\end{theorem}

We should note that in the primal subgradient algorithm   we
maintain feasibility of the coupled constraints in the problem
(\textbf{conv-DCCC}) at each iteration while for  the dual
algorithms feasibility  holds only at convergence of these
algorithms and not at the intermediate iterations. Since for control
problems the coupled constraints represent the dynamics of the
networked  system over the prediction horizon, when using a dual
algorithm these dynamics will be satisfied only at convergence. This
is a major issue when we stop at an intermediate step of a dual
based  algorithm.

There are also other dual decomposition methods  based on the
concept of augmented Lagrangians: e.g. the alternating direction
method \cite{KonLeo:96,Tse:91}, where   a quadratic penalty term
$\mu || \sum_i G_i \xb^i - g||^2$ is added to the standard
Lagrangian $L_0$. A computational drawback of this scheme is that
the quadratic penalty term is not separable in $\xb^i$. However,
this is overcome by carrying out the minimization problem in a
Gauss-Seidel fashion,
 followed by a steepest ascent update of the multipliers. In other
dual decomposition methods,  such as partial inverse method
\cite{Spi:85} or proximal point method  \cite{CheTeb:94}, for
example a term of the form $\mu \sum_i || \xb^i - \xb^i_k||^2$ is
added to the Lagrangian $L_0$. These schemes   have been shown to be
very sensitive to the value of the parameter $\mu$, with
difficulties in practice to obtain the best convergence rate. Some
heuristics for choosing $\mu$ can be found in the literature
\cite{CheTeb:94,KonLeo:96,Tse:91}. However, these heuristics have
not been formally analyzed from the viewpoint of efficiency
estimates for the general  case (linear convergence results have
been obtained e.g. only for strongly convex functions).

The new decomposition methods called     here  ``dual fast
gradient''  (\textbf{DFG})  and ``dual interior-point''
(\textbf{DIP})  obtained by smoothing the Lagrangian are more
efficient in terms of number of iterations  compared to the
classical primal or dual subgradient algorithm (see also Table  2).
We should note however, that algorithm (\textbf{DFG}) is more
appropriate than the algorithm  (\textbf{DIP}) when solving problems
where the number of coupling constraints is large, since for
(\textbf{DIP}) we need to invert at each iteration a square matrix
of dimension $n_\lambda$, where $n_\lambda$ denotes the dimension of
$\lambda$ (or equivalently the number of rows in the matrices
$G_i$).

It is also clear that the update rules in  algorithms  (\textbf{DS})
and (\textbf{DFG})  are completely distributed, according to the
communication graph between subsystems. Indeed, we recall that  the
coupling constraints $h^i(\xb^j; ~j \in \mathcal N^i) =0$ in problem
(\textbf{conv-DCCC}) are assumed to be linear,  of type $G^i
[\xb^j]_{j \in \mathcal N^i} = g_i $, i.e. we have $[G_1 \cdots G_M]
= [G^{1T} \cdots G^{MT}]^T$. Let $\lambda^i$ be the Lagrange
multipliers for the constraints $G^i [\xb^j]_{j \in \mathcal N^i} =
g_i $, and thus $\lambda = [\lambda^{1T}  \cdots \lambda^{MT}]^T$.
Then,   the main update rules in Algorithms (\textbf{DS})  and
(\textbf{DFG})  are distributed, each agent $i$ using information
only from its neighbors, e.g.:
\[ \lambda^i_{k+1}  =
  \lambda_k^i  + \alpha_k \left( G^i [\xb^j_k]_{j \in \mathcal N^i} - g_i \right). \]
However, for the algorithm  (\textbf{DIP}), the update of the
Lagrange multiplier has to be done by a central agent, i.e. in this
case we have a star-shaped topology for the communication among
subsystems. Note that for this algorithm  the sparsity of the graph
will impose sparsity  on the matrices $G_i$, which in turn will have
a strong effect on the computation of the Hessian  of the
corresponding dual function (see \cite{NecSuy:09} for more details).

\begin{table}[t]
\label{table_dcp}
\begin{center}
\begin{tabular}{|c|c|c|c|c|c|}
\hline $M$    & $N$           &   \text{nr. it. (\textbf{DS})}    &    \text{nr. it.  (\textbf{DFG})}    &    \text{nr. it.  (\textbf{DIP})} \\
\hline $10$  & $10$           &   $5.000 (0.19)$      &   $1.215 (10^{-2})$     &      $78 (10^{-4})$   \\
\hline $10$  & $20$           &    $5.000 (0.47)$     &   $ 1.873 (10^{-2})$    &     $117 (10^-{4})$   \\
\hline $10$  & $30$           &    $5.000 (0.81)$     &   $ 2.721 (10^{-2})$    &     $165 (10^{-4})$   \\
\hline
\end{tabular}
\end{center}
\caption{Distributed control problem for a network of interconnected
linear subsystems,    Example \ref{dist-ex},  where $n_i = 5, m_i=3$
and $p_i=2$ for all $i$:  we consider $M=10$ subsystems and a
prediction horizon $N=10, 20$ and $30$.  The weighted  matrices are
taken $Q_i = I_5$ and $R_i = I_2$.  By eliminating the states we
obtain the convex quadratic program \eqref{dccc-qp}  with $\xb^i
=[w_0^{iT} \cdots w_{N-1}^{iT}  \; u_0^{iT} \cdots u_{N-1}^{iT}]^T$,
where each matrix $H_i \in \rset^{N(m_i+p_i) \times N(m_i+p_i)}$ is
positive semidefinite.  In the brackets we display the accuracy
$\epsilon$.  Clearly, the dual algorithms based on smoothing
techniques (\textbf{DFG})  and (\textbf{DIP})  work much better than
classical dual subgradient algorithm  (\textbf{DS}).}
 \vspace{-0.1cm}
\end{table}

%%%%%%%%%%%%%%%%%%%%%%%%%%%%%%%%%%%%%%%%%%

\subsection{Parallel algorithms for solving    optimization problems of type (\textbf{CCDC})}

In this section  we study  parallel algorithms  for solving
optimization problems with coupled cost but decoupled constraints in
the form (\textbf{CCDC}), that e.g.  appear in the context of
cooperative control   (see Sections \ref{subsec_cooperative} and
\ref{subsec_cooperative_coupled}). A well known parallel algorithm
in linear algebra for solving systems of linear equations  is the
Jacobi algorithm that can be also used  in the context of
optimization \cite{BerTsi:89}.   Applying Jacobi algorithm, we
decompose our optimization problem of type (\textbf{CCDC}) into $M$
optimization subproblems of lower dimension.  In this algorithm each
agent updates its variable $\xb^i$ by solving a low dimension
optimization problem where the values of the rest of variables are
calculated at the previous iteration.  An extension of the Jacobi
algorithm is the Gauss-Seidel algorithm, where at each  iteration
each agent updates its variable by solving an optimization problem
for which the rest of the variables are replaced with the most
recent values computed.
\begin{algorithm}
\textbf{Algorithm Jacobi}
 \[
\xb^{i}_{k+1}= \arg \min_{\xb^i  \in \textbf{X}^i}  f(\xb^{1}_{k},
\cdots, \xb^{i-1}_{k}, \xb^{i}, \xb^{i+1}_{k}, \cdots, \xb^{M}_{k})
\]
\textbf{Algorithm Gauss-Seidel}
\[
\xb^{i}_{k+1}= \arg \min_{\xb^i \in  \textbf{X}^i}  f(\xb^{1}_{k+1},
\cdots, \xb^{i-1}_{k+1}, \xb^{i}, \xb^{i+1}_{k}, \cdots,
\xb^{M}_{k})
\]
\end{algorithm}
It is clear that in the Jacobi algorithm the  optimization
subproblems can be solved in parallel at each iteration.    The
Gauss-Seidel algorithm can be also parallelized, providing that  a
coloring scheme can be applied (see \cite{BerTsi:89} for more
details).

The convergence of these two algorithms     can be  established
under suitable contraction assumptions on the mapping $\xb - \beta
\Delta f(\xb)$  with respect to the block-maximum norm $\| \xb \| =
\max_i \| \xb^i \| / \zeta_i$ , where the $\zeta_i$'s are positive
scalars and $\xb =[ \xb^{1T} \cdots \xb^{MT} ]^T$.
\begin{theorem}  \cite{BerTsi:89}.
\label{th-jacobi} For the optimization problem (\textbf{CCDC})  we
assume  that the objective function $f$  is differentiable and
suppose that the mapping $\xb - \beta \Delta f(\xb)$  is a
contraction for some positive scalar $\beta$.  Then,  the Jacobi and
Gauss-Seidel algorithms are well defined and the sequence $\{
\xb_{k}  \}_k$ converges to the minimum of  (\textbf{CCDC}) linearly
for both iterations.
\end{theorem}

For the Gauss-Seidel algorithm, the assumptions for convergence
given in Theorem \ref{th-jacobi} can be relaxed, in particular the
contraction assumption can be replaced with a convexity assumption
on the objective function ($f$ needs to be differentiable and convex
and, furthermore, the function $f$ needs to be strictly convex
function of  $\xb^i$ when the values of  all the other components of
$\xb$ are held constant, for each $i$), see \cite{BerTsi:89} for
more details. If $f$ is not differentiable, the Jacobi or
Gauss-Seidel algorithm can fail to converge to the minimum of
(\textbf{CCDC}) because  it can stop at a non-optimal  ``corner''
point at which $f$ is non-differentiable and from which $f$ cannot
be reduced along any coordinate. The contraction assumption on the
functions $f$   for convergence of these two algorithms is usually
satisfied in many applications: see e.g. the cooperative control
problem for satellite formation   discussed in Example \ref{ex_coop}
which leads to the convex quadratic program \eqref{ccdc-qp} for
which the Hessian satisfies the contraction assumption or the
application from Section \ref{subsec_cooperative_coupled}.

In \cite{Nes:10} the optimization problem (\textbf{CCDC}) has been
solved using a coordinate descent method. The iteration $k+1$ of the
algorithm has the following form:
\begin{align*}
& \xb_{k+1}^{i_k} = \arg \min_{\xb^{i_k} \in \textbf{X}^{i_k}}
\nabla_{i_k} f(\xb_k)^T ( \xb^{i_k} - \xb_k^{i_k})
 + \frac{L_{i_k}}{2} \| \xb^{i_k} - \xb_k^{i_k} \|^2,   \\
& \xb_{k+1}^{j} = \xb_k^{j}, ~~~~~\forall j \neq i_k,
\end{align*}
where $i_k$ is chosen randomly based on a  uniform distribution.
Moreover, we assume componentwise Lipschitz continuity of the
gradient of $f$ with the Lipschitz constant $L_{i}$, for all $i=1,
\cdots, M$. In \cite{Nes:10} Nesterov proves ${\mathcal
O}(\frac{1}{\epsilon})$ rate of convergence in probability for the
coordinate descent algorithm.

For cooperative  control problems of dynamically coupled systems
(see Section \ref{subsec_cooperative_coupled}), which also leads to
optimization problems of the form (\textbf{CCDC}), various versions
of Jacobi-based algorithms have been proposed in the literature. For
example in \cite{RawSte:08,VenRaw:05,Ven:06} the authors have
proposed an algorithm of the following form:
\begin{align*}
&\overline{\xb}^i_k=\arg \min_{\xb^i \in \textbf{X}^i}
f(\xb^{1}_{k}, \cdots, \xb^{i-1}_{k}, \xb^{i}, \xb^{i+1}_{k},
\cdots,
\xb^{M}_{k}),  \\
&\xb^i_{k+1}=\alpha_i \overline{\xb}^i_k + (1-\alpha_i) \xb_k^i,
\end{align*}
where $\alpha_i$ are positive weights, summing  to $1$. In
\cite{RawSte:08,VenRaw:05,Ven:06} the authors have shown that all
the limit points of the sequence  generated by the previous
algorithm are optimal.

In \cite{CamOli:09} the authors have proposed  a decomposition of
the problem (\textbf{CCDC}) into a set of local subproblems that are
solved iteratively by a network of agents. Each subproblem ˆ is
obtained from ˆ(\textbf{CCDC})  discarding from the objective $f$
the terms that do not depend on $\xb^i$ and with the constraint set
$\textbf{X}^i$. A distributed algorithm based on the method of
feasible directions has been  proposed to generate the iterations of
the agents:
\[  \xb^i_{k+1}= \xb^i_{k} + \alpha_k^i (\hat{\xb}^i_k - \xb^i_{k}),  \]
where the local descent direction is $d_k^i = \hat{\xb}^i_k -
\xb^i_{k}$, for ˆ$ \hat{\xb}^i_k \in \textbf{X}^i$, and the step
size $\alpha_k^i$ satisfies the Armijo rule \cite{Nes:04}.  The
local iterations require relatively low effort and arrive at a
solution of (\textbf{CCDC}) at the expense of slower convergence and
high communication among neighboring agents.

\begin{table}[t]
\label{table_coo}
\begin{center}
\begin{tabular}{|c|c|c|c|c|}
\hline $M$    & $N$    &   $\sigma$      & \text{nr. it. Jacobi}         &   \text{nr. it. Gauss-Seidel}  \\
\hline $10$  & $40$   &    $0.1$           & $12.435$                       &   $3.834$     \\
\hline $10$  & $40$   &   $1$               & $1.413$                          &    $365$     \\
\hline $10$  & $40$   &   $10$             &  $174$                            &    $68$   \\
\hline
\end{tabular}
\end{center}
\caption{Cooperative control problem for satellite formation
Example \ref{ex_coop}:  we consider $M=10$ satellites and a
prediction horizon $N=40$.  The weighted  matrices are taken $Q_i =
I_3$ and $R_i = \sigma I_3$    and the accuracy of the solution
$\epsilon=10^{-3}$.  By eliminating the states we obtain the convex
quadratic program \eqref{ccdc-qp}  with $\xb^i =[\ u_0^{iT} \cdots
u_{N-1}^{iT}]^T$ and  a strongly convex objective function  having
the convexity parameter $\sigma$.  Clearly, for large $\sigma$ both
algorithms work better.}
 \vspace{-0.1cm}
\end{table}

From the Tables 1, 2 and 3 we can observe that, in order to get an
optimal solution,  we need to perform a large number of iterations.
Note however that in practical applications from control it is not
always necessary to get an optimal solution, but we can also use a
suboptimal solution that can still preserve some fundamental
properties for the system such as robustness, stability, etc.
Whenever a suboptimal solution is satisfactory we can stop the
optimization algorithm at an intermediate iteration. Note that there
exist many control strategies based on this principle of
suboptimality (see e.g. \cite{PanWri:09,ScoMay:99,VenRaw:05}).

%%%%%%%%%%%%%%%%%%%%%%%%%%%%%%%%%%%%%%%%%%

\section{Conclusions}
This paper has presented  three applications from estimation  and
process  control  for networked systems  that lead  to coupled
optimization problems with particular structure that can be
exploited in decomposition algorithms.  A  systematic framework is
then developed in the paper to explore several  parallel and
distributed  algorithms for solving such structured  optimization
problems,  each with a different tradeoff among convergence speed,
message passing amount, and distributed computation architecture.
For each application, numerical experiments on  several parallel and
distributed algorithms are provided.

%Although the theory for distributed estimation and control in complex processes is
%evolving quickly, much work remains to be done. This includes
%better tools for analyzing protocol dynamics and guaranteeing
%stable and efficient estimation and control methods, as well as better support
%for analyzing the dependencies that are introduced between the subsystems that compose the network.

%\bibliographystyle{elsarticle-harv}
%\bibliography{kuleuven_bibtexfile}

\end{document}